\documentclass[11pt]{article}
\usepackage[english]{babel}
\usepackage[cp1251]{inputenc}
\usepackage{amssymb,amsmath,amsthm,amsfonts}
\usepackage{graphicx}

\usepackage{eucal}
\usepackage{hyperref,amsthm}
\usepackage{mathptmx,bm}

\newcommand{\msc}[1]{\vspace{2mm} {\textsf{2010 MSC:} #1}}
\newcommand{\keywords}[1]{\vspace{2mm} {\textsf{Keywords}: \it #1}}

\makeatletter

\newtheorem{theorem}{Theorem}
\newtheorem{lemma}{Lemma}

\newtheorem{corollary}{Corollary}
\theoremstyle{definition}

\newtheorem{remark}{Remark}
\numberwithin{equation}{section}

\newcommand{\titre}[1]{\begin{center}
                      {\large\bf #1}
                      \end{center}}

\newcommand{\autheur}[1]{\begin{center}
                      {\bf #1}
                      \end{center}}

\newcommand{\adresse}[2]{\vspace{-10mm}\begin{center}
              {\it #1}\\{\textsf{e-mail}: #2}
                          \end{center}}

\newcommand{\dedicatory}[3]{\vspace{2mm}\begin{center}
              {\sf #1}
                          \end{center}\vspace{4mm}}

\begin{document}

\titre{\LARGE On explicit form of the Kolmogorov constant in the theory of Galton-Watson Branching Processes}
\autheur{Azam~A.~Imomov} \vspace{2mm}
    \adresse{Karshi State University; V.I.Romanovskiy Institute of Mathematics, Uzbekistan.}{imomov{\_}\,azam@mail.ru}
\autheur{Misliddin~Murtazaev} \vspace{2mm}
    \adresse{V.I.Romanovskiy Institute of Mathematics, Uzbekistan.}{misliddin1991@mail.ru}

\dedicatory{Dedicated to the fond memory of professor I.~S.~Badalbaev}

\begin{abstract}
    The paper considers the well-known Galton-Watson stochastic branching process. We are dealing with a non-critical case.
    In the subcritical case, when the mean of the direct descendants of one particle per generation of the time step is less
    than $1$, the population mean of the number of particles on the positive trajectories of the process stabilizes and
    approaches ${1}\big/\mathcal{K}$, where $ \mathcal{K}$ is the so-called Kolmogorov constant. The paper is devoted to the
    search for an explicit expression of this constant depending on the structural parameters of the process. Our reasoning is
    essentially based on the Basic Lemma, which describes the asymptotic expansion of the generating function of the distribution
    of the number of particles. An important role is also played by the asymptotic properties of the transition probabilities
    of the so-called Q-process and their property convergence to invariant measures.


\msc{60J80, 60J85}
\end{abstract}

\makeatletter
\renewcommand{\@oddhead}{\vbox{\hfill
{\textsf {\small  A.~Imomov \& M.~Murtazaev
    \qquad    On the Kolmogorov constant in Galton-Watson processes theory}}\hfill \thepage \hrule}}
\makeatother

\vspace{2mm}

\section{Introduction and main result }

    The Galton-Watson Branching Process (GWP) is a well-known classical model of population growth.
    This process describes the evolution of the population size in a system of monotype particles capable
    of death and transformation into a random number of particles of the same type. Although GWP has been
    well studied, it seems useful to discuss and clarify the well-known classical facts of GWP theory
    in more detail. In this report, we are dealing with a well-known theorem related
    to the name of A.Kolmogorov~{\cite{Kolmog38}}.

    Let the random function $Z(n)$ denote the successive population progeny in GWP at the moment
    $n\in{\mathbb{N}}_0$, where ${\mathbb{N}}_0:=\{0\}\cup\left\{{{\mathbb{N}}:=1,2,\ldots}\right\}$.
    The sequence of states $\left\{Z(n), n\in{\mathbb{N}}_0\right\}$ can be expressed as the
    following random sum of random variables:
\begin{equation*}
    Z(n+1)=\xi_{1}(n)+\xi_{2}(n)+ \cdots +\xi_{Z(n)}(n),
\end{equation*}
    where $\xi_{k}(n)$, $n,k \in {\mathbb{N}}_0$ are independent and identically distributed random variables with the
    common offspring law $\textsf{P}\left\{{\xi_{1}(1)=k}\right\}=p_k$. They are interpreted as the number of descendants
    of the $k$th particle in $n$th generation; see~{\cite[pp.~11--14]{Sev71}}. By virtue of our assumption, the stochastic
    system $\left\{Z(n)\right\}$ forms a reducible, homogeneous and discrete-time Markov chain with the state space consisting
    of two classes: $\mathcal{S}_{0}=\left\{0\right\}\cup{\mathcal{S}}$, where $\mathcal{S}\subset{\mathbb{N}}$, therein
    the state $\left\{0\right\}$ is an absorbing state, and $\mathcal{S}$ is the class of possible essential communicating
    states. Its transition probabilities
\begin{equation}              \label{A1.1}
    P_{ij}:= \textsf{P}\bigl\{{Z({n+1})=j\bigm|{Z(n)=i} }\bigr\}
    = \sum_{k_1 + \, \cdots \, + k_i=j}{p_{k_1} p_{k_2} \, \ldots \,  p_{k_i}}
\end{equation}
    for any $i,j\in\mathcal{S}$, where $p_j=P_{1j}$ and $\sum_{j\in\mathcal{S}}{p_j}=1$. This means that our GWP is
    completely defined by specifying the offspring law $\left\{{p_k, k\in\mathcal{S}}\right\}$. Conversely, any
    chain that satisfies the property \eqref{A1.1} is a GWP with the offspring law $\left\{{p_k}\right\}$;
    see {\cite[pp.1--2]{ANey}}, {\cite[p.19]{Jagers75}}. To avoid trivial cases, in what follows we
    assume that $p_k\ne 1$ and $p_0>0$, $p_0+p_1<1$.

    Considering the transition probabilities for $n$ steps
\begin{equation*}
    P_{ij} (n):=\textsf{P}_i \bigl\{{Z(n)=j}\bigr\}=\textsf{P}\bigl\{{Z({n+k})=j\bigm| {Z(k)=i}}\bigr\}
    \quad \text{for} \quad k \in {\mathbb{N}}_0,
\end{equation*}
    in conformity with \eqref{A1.1}, we find an appropriate probability generating function (GF)
\begin{equation}              \label{A1.2}
    \textsf{E}_i s^{Z(n)}:=\sum_{j\in{\mathcal{S}}}{P_{ij}(n)s^j} = \bigl[{f_n(s)}\bigr]^i,
\end{equation}
    for $s \in[0,1)$, where the GF $f_n(s)=\textsf{E}_1 s^{Z(n)}$ is $n$-fold iteration of GF
\begin{equation*}
    f(s):=\sum_{k\in{\mathcal{S}}}{p_k s^k},
\end{equation*}
    i.e. $f_{n+1}(s)=f\left(f_n(s)\right)=f_n\left(f(s)\right)$; see.~{\cite[pp.~5--6]{Harris66}}.

    Let the series $m:=\sum_{k\in\mathcal{S}}{kp_k}$ converge. Then $m=f'(1-)$ is the average of the direct
    descendants of one particle over one generation of the time step. Using the formula \eqref{A1.2}, in
    particular, one can find $\textsf{E}_1{Z(n)}=m^{n}$. In accordance with this, three classes of GWP
    are distinguished depending on the value of the parameter $m$. The process $\left\{Z(n)\right\}$ is
    called subcritical if $m<1$, critical if $m=1$, and supercritical if $m>1$ respectively.
    It is known that the sequence of vanishing probabilities $\bigl\{P_{10}(n)\bigr\}$ of one particle at time $n$
    for all classes tends monotonously to the extinction probability of the process starting from one particle which
    we will designate $q:=\lim_{n\to\infty}P_{10}(n)$. For subcritical and critical processes $q=1$,
    while in supercritical case (see.~{\cite{Vatutin2008}})
\begin{equation*}
    q=\inf\bigl\{s \in (0,1]:\, s=f(s) \bigr\}.
\end{equation*}

    In what follows, we will consider only the non-critical case, i.e., $m\neq1$, and wherever necessary,
    we will write $\textsf{E}$ and $\textsf{P}$ instead of $\textsf{E}_1$ and $\textsf{P}_1$ respectively.
    In the case under consideration $\lim_{n\to\infty}f_n(s)= q$ for all $s\in[0,1)$, and this convergence
    is uniform in $s\in[0,r]$ for any fixed $r<1$; see~{\cite[p.~53]{Sev71}}.

    Let us introduce the function $R_n(s):=q-f_n(s)$.

    In 1938 A. Kolmogorov~{\cite{Kolmog38}} established that the survival probability
    $Q(n):=\textsf{P}\bigl\{{Z(n)>0}\bigr\}=R_n(0)$ of subcritical
    process admits the asymptotic representation
\begin{equation}              \label{A1.3}
    Q(n) = \mathcal{K} m^{n} \bigl({1+o(1)}\bigr) \quad \text{as} \quad n\to\infty,
\end{equation}
    if and only if $f''(1-)<\infty$, where $\mathcal{K}$~--~a finite positive constant, called the Kolmogorov
    constant. Later, A. Nagaev and I. Badalbaev~{\cite{NBadalbaev}} improved Kolmogorov's result by proving
    the validity of the asymptotic representation~\eqref{A1.3} under a lot more weaker condition
\begin{equation*}
    \textsf{E}Z(1)\ln^{+}{Z(1)}=\sum_{k\in{\mathcal{S}}}{p_{k}k\ln{k}}<\infty.        \eqno{\textsf{$[\textsf{x}\ln{\textsf{x}}]$}}
\end{equation*}
    It follows from the representation~\eqref{A1.3}, and also noted by V.Vatutin~{\cite{Vatutin2008}},
    on positive trajectories of the process, the average number of particles population stabilizes with
    increasing generations number and approaches ${1}\large/\mathcal{K}$. Indeed, \eqref{A1.3} implies that
\begin{equation*}              \label{A1.4}
    {\frac {m^{n}} {Q(n)}}= {\frac {\textsf{E}Z(n)} {\textsf{P}\bigl\{{Z(n)>0}\bigr\}}}
    =\textsf{E}\left[{{Z(n)}\bigm|{Z(n)>0}}\right] \approx {1}\big/\mathcal{K}
    \qquad \text{as} \quad n\to\infty  \raise 0.8pt\hbox{.}
\end{equation*}
    In addition to the above, due to the last relation, the constant $\mathcal{K}$ can be interpreted as a coefficient
    of asymptotic equivalence of the average $\textsf{E}{Z(n)}=m^{n}$ of the population size to the survival probability
    of the process ${\textsf{P}\bigl\{{Z(n)>0}\bigr\}}$. The absence of an explicit expression for this constant hinders 
    the completion of a number of limit theorems for subcritical processes. In particular, this was noted long ago in 1957 
    by V. Zolotarev~{\cite{Zol57}}. Thus, it is of special interest to determine the explicit expression of this 
    equivalence coefficient as a function of the numerical parameters of the GWP.

    Our aim in this report is to get the explicit form of $\mathcal{K}$. Partly certain result on this issue is already
    available due to E.Seneta~{\cite[Theorem~2(1)]{Seneta74}} in the sense that under the condition $[\textsf{x}\ln{\textsf{x}}]$
    this constant can be calculated using the limit parameter $\mu:=\sum_{k\in{\mathbb{N}}}{k\mu_{k}}$, where
    $\mu_{k}=\lim_{n\to\infty}\textsf{P}\bigl\{{Z(n)=k \bigm| {Z(n)>0}}\bigr\}$ is the limiting-invariant
    distribution for subcritical GWP. Namely, this paper proves that
\begin{equation*}
    {\frac{m^{n}}{Q(n)}} \longrightarrow {\mu} \qquad \text{as} \quad n\to\infty,
\end{equation*}
    i.e. $\mathcal{K}={1}\big/{\mu}$.
    Consider now the random variable ${\mathcal{H}}:=\min\bigl\{{n:Z(n)=0}\bigr\}$, which denotes an extinction time of
    the process $\bigl\{Z(n)\bigr\}$ with initial state $\bigl\{Z(0)=1\bigr\}$. It is obvious that the parameter $\beta:=f'(q)$
    can be interpreted as the mean of the direct descendants of one particle in the transformed
    branching process $\bigl\{Z_{q}(n)\bigr\}$, generated by the Harris-Sevastyanov transformation $f_{q}(s)=f(qs)\big/q$. Note
    also that the process $\bigl\{Z_{q}(n)\bigr\}$ is subcritical. In this notation, the above result of E.Seneta can be
    extended to the non-critical case in the following theorem; see, also, {\cite[Lemma~2.1]{Imomov19}}.

\vspace{4mm}

\noindent{{\textbf{Theorem~S.}} {\it Let $\textsf{P}\bigl\{Z(0)=1\bigr\}=1$ and $m\neq{1}$. Then the coefficient
    of asymptotic equivalence between the population mean in the process $\bigl\{Z_{q}(n)\bigr\}$ and
    the survival probability of the process $\bigl\{Z(n)\bigr\}$ slowly stabilizes, i.e.}
\begin{equation*}
    {\frac {\textsf{E}Z_{q}(n)}{\textsf{P}\bigl\{{n<{\mathcal{H}}<\infty}\bigr\}}}=
    {\frac {\beta^{n}}{R_n(0)}}=\mathcal{L}_{\beta}\bigl({\beta^n}\bigr)   \qquad \text{\textit{as}} \quad n\to\infty,
\end{equation*}
    {\it where the function $\mathcal{L}_{\beta}\left(*\right)$ slowly varies at infinity in the sense of Karamata.
    If, in addition, the condition {$[\textsf{x}\ln{\textsf{x}}]$} holds, then}
\begin{equation*}
\lim_{n\to\infty}\textsf{P}_{i}\bigl\{{Z_{n}=k \bigm| {n<{\mathcal{H}}<\infty}}\bigr\} =: \nu_{k} <\infty,
\end{equation*}
    {\it therewith $\left\{{\nu_k, k\in\mathbb{N}}\right\}$ is a limiting-invariant
    distribution for the process $\left\{Z(n)\right\}$ such that}
\begin{equation*}
    \mathcal{K}=\lim_{n\to\infty}\mathcal{L}_{\beta}\bigl({\beta^n}\bigr)
    ={\frac {q} {\sum_{k\in{\mathbb{N}}}{k\nu_{k}}}}   \raise1.5pt\hbox{.}
\end{equation*}

\vspace{4mm}

    Further discussions show that under the Kolmogorov conditions~{\cite{Kolmog38}}, the constant $\mathcal{K}$ 
    can be explicitly calculated using the structural parameters (moments) of the process $\bigl\{Z(n)\bigr\}$.
    Below we formulate the main result of the paper, in which the explicit form of $\mathcal{K}$
    is found depending on the second factorial moment $f''(1-)$.

\begin{theorem}              \label{IMTh:1}
    Let $m\ne1$ and $2b_{q}:=f''(q)<\infty$. Then
\begin{equation*}
    \mathcal{K}={\frac {q} {\,1+q\gamma\,}}\,   \raise1.5pt\hbox{,}
\end{equation*}
    where $\gamma={{b_{q}}\big/{\left(\beta-\beta^{2}\right)}}$.
\end{theorem}

\begin{corollary}              \label{IMCor:1}
    Let $m<1$. If $2b:=f''(1-)<\infty$, then
\begin{equation*}
    \mathcal{K}={\frac {1} {\,1+\gamma\,}}\,  \raise1.5pt\hbox{,}
\end{equation*}
    where $\gamma={{b}\big/{\left(m-m^{2}\right)}}$.
\end{corollary}

    Section~\ref{MySec:2} we devote to the proof of the Theorem~\ref{IMTh:1}.

\section{The proof of Theorem~\ref{IMTh:1}}     \label{MySec:2}

    We divide the proof of the theorem into several steps.

\subsection{A defective but important Lemma}     \label{MySubsec:2.1}
    The mean value theorem gives
\begin{equation}              \label{A2.1}
    R_{n+1}(s)=f'\bigl({\xi_n(s)}\bigr)R_n(s),
\end{equation}
    where $\xi_n(s)=q-\theta R_n(s)$ and $0<\theta <1$. We see that if $s\in[0,q)$ then $R_n(s)>0$, therefore
    $\xi_n(s)<q$. Since the GF $f(s)$ and its derivatives are monotonically nondecreasing, successive application 
    of \eqref{A2.1} leads to the inequality $R_n(s)<q\beta^n$. Thence
\begin{equation*}
    q-q\beta^n<\xi_n(s)<q  \qquad \text{for} \quad s\in[0,q).
\end{equation*}
    Accordingly
\begin{equation*}
    {\frac{R_{n+1}(s)} {\beta}}<R_n(s)<{\frac{R_{n+1}(s)} {f'\bigl(q-q\beta^n\bigr)}}  \raise 1.5pt\hbox{.}
\end{equation*}
    On the other hand, for all $s\in[q,1)$ we see $R_n(s)<0$, so that $\xi_n(s)=q+\theta\left|R_n(s)\right|>q$.
    And in this case, successively applying the formula \eqref{A2.1} and, taking into account the properties of 
    $f(s)$, we obtain the inequality $(q-1)\beta^n <R_n(s)$ or the same as $\left|R_n(s)\right| <(1-q)\beta^n$. Hence
\begin{equation*}
    q <\xi_n(s)<q+(1-q)\beta^n  \qquad \text{for} \quad s\in[q,1).
\end{equation*}
    Then
\begin{equation*}
    {\frac{R_{n+1}(s)} {f'\bigl(q+(1-q)\beta^n\bigr)}}<R_n(s)<{\frac{R_{n+1}(s)} {\beta}}  \raise 1.5pt\hbox{.}
\end{equation*}

    Based on recent results, we conclude that
\begin{equation}              \label{A2.2}
    {\frac{R_{n+1}(s)} {f'\bigl(q_{1}(n)\bigr)}}<R_n(s)<{\frac{R_{n+1}(s)} {f'\bigl(q_{0}(n)\bigr)}}
    \qquad \text{for all} \quad s\in[0,1),
\end{equation}
    where
\begin{equation*}
    q_{k}(n):=q+(k-q)\beta^n \qquad \text{for} \quad  k=0,\,1.
\end{equation*}

    In turn, by the Taylor formula and by iterating over $f(s)$ we have the following relation:
\begin{equation}              \label{A2.3}
    R_{n+1}(s)=\beta R_n (s)-{\frac{f''\bigl(\zeta_n(s)\bigr)}{2}}R_n^2(s)
\end{equation}
    for all all $s\in[0,1)$, herein $\zeta_n(s)$ is so that
\begin{equation*}
     q_{0}(n)<\zeta_n(s)< q_{1}(n)   \qquad \text{for all} \quad s\in[0,1).
\end{equation*}
    Accordingly, by the monotone non-decreasing property of GF we obtain
\begin{equation}              \label{A2.4}
    f''\bigl(q_{0}(n)\bigr)<f''\bigl(\zeta_n(s)\bigr)<f''\bigl(q_{1}(n)\bigr).
\end{equation}
    Combining \eqref{A2.2}--\eqref{A2.4} implies
\begin{equation*}
    {\frac{f''\bigl(q_{0}(n)\bigr)}{2f'\bigl(q_{1}(n)\bigr)}}{R_n(s)R_{n+1}(s)} < \beta R_n (s)-R_{n+1}(s)<
    {\frac{f''\bigl(q_{1}(n)\bigr)} {2f'\bigl(q_{0}(n)\bigr)}}{R_n(s)R_{n+1}(s)}.
\end{equation*}
    Multiplying these inequalities to $1\big/\bigl({R_n(s)R_{n+1}(s)}\bigr)$ gives us
\begin{equation}              \label{A2.5}
    {\frac{f''\bigl(q_{0}(n)\bigr)}{2f'\bigl(q_{1}(n)\bigr)}}<{\frac{\beta}{R_{n+1}(s)}}
    - {\frac{1}{R_n(s)}}<{\frac{f''\bigl(q_{1}(n)\bigr)} {2f'\bigl(q_{0}(n)\bigr)}}  \raise 1.5pt\hbox{.}
\end{equation}
    Repeated application inequalities \eqref{A2.5} leads to the following ones:
\begin{equation*}
    {\frac{1}{2}}\sum_{k=0}^{n-1}{\frac{f''\bigl(q_{0}(k)\bigr)}{f'\bigl(q_{1}(k)\bigr)}{\beta^{k}}}
    <{\frac{\beta^{n}}{R_{n}(s)}}-{\frac{1}{1-s}}<
    {\frac{1}{2}}\sum_{k=0}^{n-1}{\frac{f''\bigl(q_{1}(k)\bigr)}{f'\bigl(q_{0}(k)\bigr)}{\beta^{k}}}  \raise 1.5pt\hbox{.}
\end{equation*}
    Taking limit as $n \to \infty $ from here we have estimation
\begin{equation}              \label{A2.6}
    {\frac{\Delta_{1}}{2}} \leq
    \lim_{n\to\infty}\left[{\frac{\beta^{n}}{R_{n}(s)}}-{\frac{1}{1-s}}\right]
    \leq {\frac{\Delta_{2}}{2}}  \raise 1.5pt\hbox{,}
\end{equation}
    where
\begin{equation*}
    \Delta_1:= \sum_{k=0}^{\infty}{\frac{f''\bigl(q_{0}(k)\bigr)}{f'\bigl(q_{1}(k)\bigr)}{\beta^{k}}}
         \qquad {\text{and}} \qquad
    \Delta_2:= \sum_{k=0}^{\infty}{\frac{f''\bigl(q_{1}(k)\bigr)}{f'\bigl(q_{0}(k)\bigr)}{\beta^{k}}}.
\end{equation*}
    Evidently, last two series converge because of
    $0<q_{0}(n)<q<q_{1}(n)<1$ for all $n\in\mathbb{N}_{0}$. Designating
\begin{equation*}
    {\frac{1} {A_1(s)}}:= {\frac{1} {q-s}}+{\frac{\Delta_1} {2}}
        \qquad {\text{and}} \qquad
    {\frac{1} {A_2(s)}}:= {\frac{1} {q-s}}+{\frac{\Delta_2} {2}} \raise 1.5pt\hbox{,}
\end{equation*}
    we rewrite the relation \eqref{A2.6} as follows:
\begin{equation}              \label{A2.7}
    {\frac{1} {A_1(s)}} \le \lim_{n\to\infty}{\frac{\beta^n} {R_n(s)}} \le {\frac{1} {A_2(s)}}  \raise 1.5pt\hbox{.}
\end{equation}
    Clearly 
\begin{equation*}
    {\frac{1} {A_2(s)}} - {\frac{1} {A_1(s)}} = {\frac{\Delta_2 - \Delta_1} {2}} < \infty.
\end{equation*}
    In turn, we see that ${{\beta^n} \big/{R_n(s)}}$ monotonously increases for all $s\in[0,q)$
    as $n\to\infty$ and monotonously increases and decrease in kind for all $s\in[q,1)$. Therefore
    ${\mathcal{A}_{\delta}(s)}:=\lim_{n\to\infty}{{\beta^n}\big/{R_n(s)}}$ exists and in accordance with
    inequalities \eqref{A2.7} there is a positive variable $\delta\in\left[\Delta_1, \Delta_2\right]$
    such that
\begin{equation}              \label{A2.8}
    {\frac{1} {\mathcal{A}_{\delta}(s)}}:= {\frac{1} {q-s}}+{\frac{\,\delta\,} {2}}  \raise 1.5pt\hbox{.}
\end{equation}

    So we established the following statement.
\begin{lemma}               \label{AILem:1}
    If $m\ne{1}$ and $f''(q)<\infty$, then
\begin{equation}              \label{A2.9}
    {\frac{R_n(s)}{\beta^{n}}} \longrightarrow {\mathcal{A}_{\delta}(s)}
    \qquad \text{as} \quad n\to\infty,
\end{equation}
    where the function ${\mathcal{A}_{\delta}}(s)$ is defined in \eqref{A2.8}.
\end{lemma}

    The core and only defect of Lemma~\ref{AILem:1} is the lack of an explicit $\delta$ expression 
    in \eqref{A2.8}. In Section~\ref{MySubsec:2.2} we will eliminate this defect.

\subsection{The Q-process contribution}     \label{MySubsec:2.2}

    We begin by recalling the so-called Q-process which is an irreducible homogeneous-discrete-time Markov
    chain $\bigl\{W(n)\bigr\}$ ith the state space $\mathcal{E}\subset{\mathbb{N}}$. The transition probabilities
    of Q-process are 
\begin{equation}              \label{A2.10}
    {\mathcal{Q}}_{ij}(n):=\textsf{P}\left\{{W({n+k})=j}\bigm|{W(k)=i}\right\}={\frac{jq^{j-i}}{i\beta^n}}P_{ij}(n)
    \qquad \text{for all} \quad i,j \in {\mathcal{E}},
\end{equation}
    and for any $n,k \in {\mathbb{N}}$; see.~\cite[Sec.~I, \S14]{ANey}.
    Put into consideration a GF
\begin{equation*}
    w_n^{(i)}(s):=\sum\limits_{j\in{\mathcal{E}}} {{\mathcal{Q}}_{ij}(n)s^j}.
\end{equation*}
    Then from \eqref{A1.2} and \eqref{A2.10} we have
\begin{eqnarray}
    w_n^{(i)}(s) \nonumber
    & = & \sum\limits_{j \in {\mathcal{E}}} {{\frac{jq^{j-i}}{i\beta^n}}P_{ij}(n)s^j} \\
    & = & {\frac{q^{1-i}s}{i\beta^n}}\sum\limits_{j\in{\mathcal{E}}}{P_{ij}(n)(qs)^{j-1}}
    = {\frac{qs}{i\beta^n}}{\frac{\partial}{\partial{x}}}\left[{\left({{\frac{f_n(x)}{q}}} \right)^i}\right]_{x=qs}. \nonumber
\end{eqnarray}
    Last formula is convenient for using in a following form:
\begin{equation}              \label{A2.11}
    w_n^{(i)}(s)=\left[{{\frac{f_n(qs)} {q}}} \right]^{i-1}w_n(s),
\end{equation}
    where the GF $w_n(s):=w_n^{(1)}(s)=\textsf{E}\left[{s^{W(n)} \bigm| {W(0) =1}} \right]$ has a form of:
\begin{equation}              \label{A2.12}
    w_n (s) = s{\frac{f'_n(qs)} {\beta^n}}
    \qquad \text{for all} \quad n \in  {\mathbb{N}}.
\end{equation}
    Since $f_n(s)\to q$ uniformly in $s\in[0,r]$ for any fixed $r<1$ as ${n\to\infty}$, it follows from
    \eqref{A2.11} and \eqref{A2.12} that ${{\mathcal{Q}}_{ij}(n)} \big/ {{\mathcal{Q}}_{1j}(n)}\to 1$
    as infinitely growth the number of generations.

    Application of iteration for $f(s)$ in the relation \eqref{A2.11} leads us to the following functional equation:
\begin{equation}              \label{A2.13}
    w_{n+1}^{(i)}(s)={\frac{w(s)} {f_{q}(s)}}w_n^{(i)}\bigl({f_{q}(s)}\bigr),
\end{equation}
    where $w(s):=w_1(s)$ and $f_{q}(s)= {f(qs)\big/ q}$. Thus, Q-process is completely defined by setting the GF
\begin{equation}              \label{A2.14}
    w(s) = s{\frac{f'(qs)} {\beta}}\raise 1.5pt\hbox{.}
\end{equation}
    An evolution of the Q-process is in essentially regulated by the structural parameter $\beta >0$. In fact,
    as it has been shown in \cite[p.~59,~Theorem~2]{ANey}, that if $\beta <1$ then $\mathcal{E}$ is positive
    recurrent and, $\mathcal{E}$ is transient if $\beta =1$. On the other hand, it is easy to be convinced
    that positive recurrent case $\beta <1$ of Q-process corresponds to the non-critical case $m\neq{1}$
    of GWP. Note that $\beta \leq{1}$ and nothing but.

    Assume that $\alpha :=w'(1-)<\infty$ in the case $\beta <1$. Then differentiating \eqref{A2.14}
    on the point $s=1$ we obtain $\alpha = 1+\left({1-\beta}\right)\gamma_{q}$, where
\begin{equation*}
    \gamma_{q}:={\frac{qf''(q)}{\left({\beta-\beta^2}\right)}} \raise 1.5pt\hbox{.}
\end{equation*}
    Further, it follows from \eqref{A2.11} and \eqref{A2.12} that
    ${\textsf{E}}_i W(n)=\left({i-1} \right)\beta^n+\textsf{E}W(n)$, where
\begin{equation*}
    \textsf{E}W(n)= 1+\gamma_{q}\cdot\bigl({1-\beta^n}\bigr).
\end{equation*}
    
    It is known \cite[p.~59,~Theorem~2(iv)]{ANey} that in this case there exists an invariant measure 
    $\bigl\{\pi_j\bigr\}$ with respect to the probabilities ${{\mathcal{Q} }_{ij}(n)}$ such that
\begin{equation}              \label{A2.15}
    \pi_j:=\lim_{n\to\infty}{\mathcal{Q}}_{ij}(n)=jq^{j-1}\nu_{j}  \qquad \text{for all} \quad  i,j\in{\mathcal{E}},
\end{equation}
    where $\{\nu_{j}\}$ are coefficients in a power series expansion of the limit GF
\begin{equation*}
    Q(s):=\lim_{n\to\infty}{\frac{f_n(s)-q}{\beta^n}};
\end{equation*}
    see~\cite[p.~41,~Theorem~3]{ANey}.
    In conformity with our designation and by Lemma~\ref{AILem:1} we see $Q(s)=-\mathcal{A}_{\delta}(s)$.
    Then $\mathcal{A}_{\delta}(s)=-\sum_{j\in{\mathcal{E}}}{\nu_js^j}$. Thus, interpreting 
    the statement \eqref{A2.15} in the context of GF, we conclude that there exists 
    a limit GF $\pi(s):=\sum_{j\in{\mathcal{E}}}{\pi_j s^j}$ such that
\begin{equation}              \label{A2.16}
    \pi(s)=\lim_{n\to \infty}w_n^{(i)}(s)=-s\mathcal{A}'_{\delta}(qs).
\end{equation}
    for all $s\in [0,1)$. On the other hand, taking limit as $n\to\infty$ in equation~\eqref{A2.13} with a combination of
    equations~\eqref{A2.11} and \eqref{A2.12}, leads us to the following Schr\"{o}der type functional equation:
\begin{equation}              \label{A2.17}
    \pi(s)={\frac{w(s)} {f_{q}(s)}}\pi\bigl({f_{q}(s)}\bigr).
\end{equation}
    The equation~\eqref{A2.17} entails
    $\pi_j=\sum_{i\in{\mathcal{E}}}{\pi_i{\mathcal{Q}}_{ij}(n)}$ for $n\in{\mathbb{N}}$ and $j\in{\mathcal{E}}$.

    Now, due to the form of \eqref{A2.8} and from the relation \eqref{A2.16} immediately follows, that
    $\pi(1)=1$. The last argument is equivalent to that $\bigl\{\pi_j\bigr\}$ represents an invariant
    distribution. Simultaneously differentiating the equation \eqref{A2.17} and taking $s=1$ we obtain
    the mean of distribution $\bigl\{\pi_j\bigr\}$ as 
\begin{equation*} 
    \pi'(1-)=1+\gamma_{q}. 
\end{equation*}
    At the same time, the relation \eqref{A2.16} implies that $\pi'(1-)=1+q\delta$. Hence
\begin{equation}              \label{A2.18}
    \delta={\frac{\gamma_{q}} {\,q\,}} = {\frac{f''(q)} {\beta\bigl(1-\beta\bigr)}}   \raise 1.5pt\hbox{.}
\end{equation}

\subsection{Basic Lemma}     \label{MySubsec:2.3}

    From the representation \eqref{A2.8} and the equality \eqref{A2.18} one can finally obtain an explicit expression for 
    the limit function $\lim_{n\to\infty}{{R_n(s)}\big/{\beta^n}}$ depending on $\beta$ and ${f''(q)}$ for all $s\in [0,1)$. 
    Thus, we have proved the following Basic Lemma, in which the deficiency of the lemma~\ref{AILem:1} is eliminated.
\begin{lemma}               \label{AILem:2}
    If $m\ne{1}$ and $2b_{q}:=f''(q)<\infty$, then
\begin{equation*}
    R_n(s)={\mathcal{A}_{\gamma}(s)} \beta^n\bigl({1+o(1)}\bigr)
    \qquad \text{as} \quad n\to\infty,
\end{equation*}
    where
\begin{equation}              \label{A2.19}
    {\frac{1} {\mathcal{A}_{\gamma}(s)}}= {\frac{1} {q-s}}+{\gamma},
\end{equation}
    and $\gamma={{b_q}\big/{\left(\beta-\beta^{2}\right)}}$.
\end{lemma}

    Finally, the statement of Theorem~\ref{IMTh:1} follows from Lemma~\ref{AILem:2}, taking in $s=0$:
\begin{equation*}
    0<\mathcal{K}={\mathcal{A}_{\gamma}(0)}={\frac{q}{1+q\gamma}}<\infty.
\end{equation*}
    Hence, taking $q=1$, we obtain the assertion of Corollary~\ref{IMCor:1}.

\section {Attendant remarks}

\begin{remark}               \label{AIRem:1}
    Due to the expression in \eqref{A2.19}, we obtain the following properties of the function ${\mathcal{A}_{\gamma}(s)}$:
\begin{itemize}
\item[$\blacktriangleright$]  ${\mathcal{A}_{\gamma}(q)}=0$;

\item[$\blacktriangleright$]  ${\mathcal{A}'_{\gamma}(q)}=-1$;

\item[$\blacktriangleright$]  it asymptotically satisfies to the Schr\"{o}der functional equation, i.e.
\begin{equation*}
    {\mathcal{A}_{\gamma}\bigl({f_n(qs)}\bigr)}=\beta^n {\mathcal{A}_{\gamma}(qs)}\bigl({1+o(1)}\bigr)
    \qquad \text{as} \quad n\to\infty
\end{equation*}
    for all $s\in[0,1)$.
\end{itemize}
    These properties are in full compliance with the properties established in \cite[Sec.~I,~\S11]{ANey} and
    also in {\cite{Imomov14a}} for a continuous-time Markov branching process; see also {\cite{Imomov12}}.
\end{remark}

\begin{remark}               \label{AIRem:2}
    The following differential analogue of the Basic Lemma plays more important role
    in the theory of non-critical GWP, which we will write out from last findings:
    \textit{if $m\neq1$ and $f''(q)<\infty$, then
\begin{equation}              \label{A3.1}
    {\frac{\partial{R_n(s)}}{\partial{s}}} =
    -{\frac{\mathcal{A}^{2}_{\gamma}(s)}{(q-s)^{2}}} \, \beta^n \cdot \bigl({1+o(1)}\bigr)
    \qquad \text{as} \quad n\to\infty,
\end{equation}
    where ${\mathcal{A}_{\gamma}(s)}$ is defined in \eqref{A2.19}.}
    Applying \eqref{A3.1} directly, setting $s=0$ there, we obtain the following asymptotic expansion:
\begin{equation*}
    \beta^{-n}P_{11}(n) = {\frac{1}{q^{2}}}\,{\mathcal{K}^{2}} \cdot \bigl({1+o(1)}\bigr)
    \qquad \text{as} \quad n\to\infty.
\end{equation*}
\end{remark}

\begin{remark}               \label{AIRem:3}
    Apparently, arguments like the last one will allow one to calculate $\mathcal{K}$ for a continuous-time
    Markov branching process in which the second factorial moment of the branching rate law is finite.
\end{remark}

\end{document}